\documentclass[a4paper,12pt]{article}
\usepackage[margin=1in]{geometry}
\usepackage{amsmath,amsthm,amssymb,amsfonts}
\usepackage{hyperref}

\newtheorem{proposition}{Proposition}
\newtheorem{theorem}{Theorem}
\newtheorem{lemma}{Lemma}

\title{A sharper log-convexity inequality for Bell numbers}
\author{Vuong Bui\thanks{\texttt{bui.vuong@yandex.ru}}}
\date{}

\begin{document}

\maketitle

\begin{abstract}
    We prove a stronger version of the log-convexity inequality for the Bell numbers $B_n$. In particular, for $n\ge 5$, we have
    \[
        B_{n+1}B_{n-1} - (B_n)^2 \ge \sum_{i=1}^{n} F_i (B_{n-i})^2,  
    \]
    where $F_i$ is the $i$-th Fibonacci number with $F_0=F_1=1$.
    The simple proof is mostly combinatorial with elementary inequalities.
\end{abstract}

\section{Introduction}
Let $B_n$ denote the Bell number, that is the number of partitions of $[n]=\{1,2,\dots,n\}$, with the convention that $B_0=1$. More information can be found at the entry OEIS sequence A000110~\cite{OEIS_A000110}. It is well known that $B_n$ is log-convex, that is $B_{n+1}B_{n-1}\ge (B_n)^2$ for every $n\ge 1$. Engel \cite{engel1994average} gave the first proof and Canfield later \cite{canfield1995engel} gave a simpler one. A refinement of Engel's inequality was later obtained by Alzer \cite{alzer2019engel}, who showed that $B_{n-1}B_{n+1}-B_n^2$ can be represented as an infinite series with nonnegative terms. Our result gives a different type of refinement, namely an explicit finite lower bound for this log-convexity gap in terms of smaller Bell numbers.

\begin{theorem} \label{thm:main}
    For $n\ge 5$, we have
    \[
        B_{n+1}B_{n-1}\ge \sum_{i=0}^{n} F_{i} (B_{n-i})^2,  
    \]
    where $F_i$ is the $i$-th Fibonacci number with $F_0=F_1=1$.
\end{theorem}
Note that the inequality fails for $n=2,3,4$. However, if we drop some terms or shift some coefficients, we can obtain a version that works for every $n\ge 1$.

We first prove a weaker version with only the first two terms on the right hand side, but still stronger than the log-convexity $B_{n+1}B_{n-1}\ge (B_n)^2$. It is to show how one can easily prove the log-convexity with combinatorial arguments and elementary inequalities.

\begin{proposition} \label{prop:simple}
    For every $n\ge 1$, we have
    \[
        B_{n+1}B_{n-1}\ge (B_n)^2 + (B_{n-1})^2.
    \]
\end{proposition}

\begin{proof}
    For each partition $p$ of $[n]$, we call the partition $s$ after removing $n$ the \emph{base} of $p$. Since the block containing $n$ can be precisely $\{n\}$, we add an empty set to each of the base in case there is not. For example, for $n=5$, the base of $\{\{1,2,3\},\{4,5\}\}$ is $\{\{1,2,3\},\{4\}, \emptyset\}$ and the base of $\{\{1,2,3\},\{4\}, \{5\}\}$ is also $\{\{1,2,3\},\{4\}, \emptyset\}$. One can see that the number of possible bases is precisely the number of partitions of $[n-1]$.
    Let $B_n^{(s)}$ denote the number of partitions of $[n]$ with base $s$. By the Cauchy--Schwarz inequality,
    \[
        (B_n)^2 = \left(\sum_s B_n^{(s)}\right)^2 \le \left(\sum_s 1\right)\left(\sum_s (B_n^{(s)})^2\right) = B_{n-1} \left(\sum_s (B_n^{(s)})^2\right).
    \]
    For each pair of partitions of $[n]$ with the same base $s$, we construct a partition of $[n+1]$ as follows. Rename the element $n$ in the second partition to $n+1$, and then merge the two copies of the common base $s$, keeping track of the blocks to which $n$ and $n+1$ are attached. If these two blocks coincide, then $n$ and $n+1$ are put in the same block. This construction is injective. The partitions not obtained in this way are precisely those containing both singleton blocks $\{n\}$ and $\{n+1\}$, and there are $B_{n-1}$ such partitions.
    Totally,
    \[
        (B_n)^2\le B_{n-1} (B_{n+1} - B_{n-1}).
    \]
    The conclusion follows.
\end{proof}

We discuss Canfield's proof in Section \ref{sec:canfield}. We will then give a weaker version of Theorem \ref{thm:main} in Section \ref{sec:stronger} for explaining the ideas better. Theorem \ref{thm:main} is proved in Section \ref{sec:strongest}.

\section{Some discussions on Canfield's proof}
\label{sec:canfield}
Let us revisit the simplest proof of the log-convexity in literature, which was given by Canfield in \cite{canfield1995engel}.
In the proof, we let $S_n$ denote the number of blocks over all partitions of $[n]$.
A partition of $[n]$ can be extended to a partition of $[n+1]$ by either (i) adding the block $\{n+1\}$ to the partition or (ii) inserting $n+1$ into one of the pre-existing blocks. Applying to all partitions of $[n]$, the resulting partitions cover all the partitions of $[n+1]$. The number of the resulting partitions by following (i) is $B_n$ while the number of those by following (ii) is $S_n$. Therefore, the log-convexity of $B_n$, which can be rewritten as $B_{n+1}/B_n\ge B_n/B_{n-1}$, is equivalent to
\begin{equation} \label{eq:mono}
    \frac{S_n}{B_n}\ge \frac{S_{n-1}}{B_{n-1}},
\end{equation}
since $B_{n+1}/B_n = (B_n+S_n)/B_n = 1+S_n/B_n$.
After that, it is argued that $S_n/B_n$, which is the average number of blocks per partition, should be increasing due to the fact that the more blocks in a partition of $[n]$, the more resulting partitions of $[n+1]$, with as many as or one more blocks in each partition.
This observation is so intuitive that Canfield did not provide a detailed verification. It turns out that by clarifying it, we can obtain a stronger result, identical to the one in Proposition \ref{prop:simple}. As one can guess, some application of the Cauchy--Schwarz inequality is involved.

Indeed, let $\mathcal S_n$ denote the multiset of the number of blocks in the partitions of $[n]$. In other words, $S_n$ is the sum of all the elements in $\mathcal S_n$. Let $Q=\mathcal S_{n-1}$, the monotonicity of $S_n/B_n$ in \eqref{eq:mono} is equivalent to
\[
    \frac{\sum_{k\in Q} k^2 +k + 1}{\sum_{k\in Q} k + 1} \ge \frac{\sum_{k\in Q} k}{\sum_{k\in Q} 1}.
\]
We explain why $S_n = \sum_{k\in Q} k^2 +k + 1$: for each partition of $[n-1]$ with $k$ blocks, we extend it to $k$ partitions, each with the same number of $k$ blocks, and another partition with $k+1$ blocks. The reason for $B_n=\sum_{k\in Q} k + 1$ is similar, when we do not care about the number of blocks but just the number of partitions.

Rewriting the inequality with products instead fractions, it becomes
\[
    \left(\sum_{k\in Q} k\right)\left(\sum_{k\in Q} k + 1\right) \le \left({\sum_{k\in Q} 1}\right)\left(\sum_{k\in Q} k^2 +k + 1\right).
\]
This inequality just follows from the Cauchy--Schwarz inequality, hence \eqref{eq:mono} follows.

One can notice that we have not used the full strength of the Cauchy--Schwarz inequality. Since $k(k+1) = (k^2+k+1)-1$, we can also obtain
\[
    \left(\sum_{k\in Q} k\right)\left(\sum_{k\in Q} k + 1\right) \le \left({\sum_{k\in Q} 1}\right)\left(\sum_{k\in Q} k^2 +k + 1\right) - \left({\sum_{k\in Q} 1}\right)\left(\sum_{k\in Q} 1\right),
\]
which is equivalent to
\[
    B_n S_{n-1} \le B_{n-1} S_n - (B_{n-1})^2.
\]
In other words,
\[
    \frac{S_{n-1}}{B_{n-1}} + \frac{B_{n-1}}{B_n} \le \frac{S_n}{B_n}.
\]
Adding one to both sides yields
\[
    \frac{B_n}{B_{n-1}} + \frac{B_{n-1}}{B_n} \le \frac{B_{n+1}}{B_n},
\]
which is the same conclusion as Proposition \ref{prop:simple}.

\section{Uniform coefficients}
\label{sec:stronger}
We now prove a stronger conclusion than Proposition \ref{prop:simple}.
\begin{proposition} \label{prop:stronger}
    For every $n\ge 3$, we have
    \[
        B_{n+1}B_{n-1} \ge  (B_n)^2 + B_nB_{n-2}.
    \]
\end{proposition}
Note that it is stronger due to $(B_{n-1})^2 \le B_nB_{n-2}$. However, it works only for $n\ge 3$. Another way to see it as a stronger result is via the following Theorem \ref{thm:uniform}, which we will derive as a corollary of Proposition \ref{prop:stronger} later.
\begin{theorem} \label{thm:uniform}
    For every $n\ge 1$, we have
    \[
        B_{n+1}B_{n-1}\ge (B_n)^2 + (B_{n-1})^2 + \dots + (B_1)^2.
    \]    
\end{theorem}

In order to prove Proposition \ref{prop:stronger} (and later results), we need the following lemma.

\begin{lemma} \label{lem:recurrence}
    Let $\mathcal S_n$ denote the multiset of the number of blocks in the partitions of $[n]$.
    For every positive $n$ and every nonnegative $t$, we have
    \[
        B_{n+t} = \sum_{k\in\mathcal S_n} \sum_{s=0}^t \binom{t}{s} k^s B_{t-s}.
    \]
\end{lemma}
\begin{proof}
    Consider a partition of $[n]$ with $k$ blocks. If we extend it to a partition of $[n+t]$, we need to pick some $s$ elements of $[n+t]\setminus[n]$ to join the $k$ blocks and the remaining $t-s$ elements would occupy extra blocks. The conclusion follows from the fact that the relation is a bijection and
    \begin{itemize}
        \item there are $\binom{t}{s}$ choices for the $s$ elements,
        \item there are $k^s$ choices for the $s$ elements to join the $k$ blocks,
        \item there are $B_{t-s}$ choices for the remaining $t-s$ elements to occupy extra blocks.\qedhere
    \end{itemize}
\end{proof}

We are now ready to prove Proposition \ref{prop:stronger}.
\begin{proof}[Proof of Proposition \ref{prop:stronger}]
    Applying Lemma \ref{lem:recurrence} for $t=0,1,2,3$, we have
    \begin{align*}
        B_n &= \sum_{k\in \mathcal S_n} 1,\\
        B_{n+1} &= \sum_{k\in \mathcal S_n} k+1,\\
        B_{n+2} &= \sum_{k\in \mathcal S_n} (k+1)^2+1,\\
        B_{n+3} &= \sum_{k\in \mathcal S_n} (k+1)^3+3(k+1)+1.
    \end{align*}

    Let $Q= \mathcal S_{n-2}$. The rewritten conclusion
    \[
        B_{n+1}B_{n-1}\ge B_n(B_n + B_{n-2})
    \]
    is equivalent to
    \[
        \left(\sum_{k\in  Q} (k+1)^3+3k+3+1\right)
        \left(\sum_{k\in  Q} k+1\right)
        \ge
        \left(\sum_{k\in  Q} (k+1)^2+1\right)
        \left(\sum_{k\in  Q} (k+1)^2+1+1\right).
    \]
    
Proving the inequality is standard as follows.
Let $x=k+1$. Since $k$ is a positive integer, we have $x\ge 2$. The desired inequality becomes
\[
\left(\sum_x x^3+3x+1 \right)\left(\sum_x x\right)\ge \left(\sum_x x^2+1 \right)\left(\sum_x x^2+2 \right),
\]
where the sums are over the multiset $Q+1 = \{k+1:k\in Q\}$.

Let $D$ be the difference between the left-hand side and the right-hand side. Then
\[
D=\sum_{x,y}\left[(x^3+3x+1)y-(x^2+1)(y^2+2)\right].
\]
Doubling and symmetrizing gives
\begin{align*}
2D&=\sum_{x,y}\left[(x^3+3x+1)y+(y^3+3y+1)x-(x^2+1)(y^2+2)-(y^2+1)(x^2+2)\right] \\
&=\sum_{x,y}\left[(xy-3)(x-y)^2+x+y-4\right].
\end{align*}
Since $x,y\ge 2$, we have
\[
xy-3\ge 0\qquad\text{and}\qquad x+y-4\ge 0.
\]
Hence every term in the final sum is nonnegative, so $2D\ge 0$. Therefore $D\ge 0$. The conclusion follows.
\end{proof}

We can now prove Theorem \ref{thm:uniform}. 
\begin{proof} [Proof of Theorem \ref{thm:uniform}]
We assume $n\ge 3$, since otherwise Theorem \ref{thm:uniform} can be deduced from Proposition \ref{prop:simple} for $n=1,2$. Applying Proposition \ref{prop:stronger} iteratively from $n$ to $3$, we have
\begin{align*}
    B_{n+1}B_{n-1}&\ge (B_n)^2 + B_nB_{n-2} \\
    &\ge (B_n)^2 + (B_{n-1})^2 + B_{n-1}B_{n-3} \\
    &\ge \dots \\
    &\ge (B_n)^2 + (B_{n-1})^2 + \dots + (B_4)^2 + B_4B_2\\
    &\ge (B_n)^2 + (B_{n-1})^2 + \dots + (B_4)^2 + (B_3)^2 + B_3B_1\\
    &\ge (B_n)^2 + (B_{n-1})^2 + \dots + (B_4)^2 + (B_3)^2 + (B_2)^2 + (B_1)^2,
\end{align*}
where the last step is due to $B_3B_1 = (B_2)^2 + (B_1)^2 = 5$ (or one can apply Proposition \ref{prop:simple}, which works for every $n$).
\end{proof}

\section{Exponential coefficients}
\label{sec:strongest}
While Theorem \ref{thm:uniform} is quite stronger than the log-convexity, one can wonder if we can make the coefficients of $(B_{n-i})^2$ grow exponentially in terms of $i$. This section shows that it is possible with the Fibonacci numbers as the coefficients, which is the conclusion of Theorem \ref{thm:main}. In fact, if we extend the arguments in the way below, we can improve the results even further. However, the proof will become more complicated. Therefore, we stop at the nice coefficients of Fibonacci numbers.

We recall Theorem \ref{thm:main}: For $n\ge 5$, we have
    \[
        B_{n+1}B_{n-1}\ge (B_n)^2 + \sum_{i=1}^{n} F_{i} (B_{n-i})^2. 
    \]
    
    The key immediate step is the following proposition, which can be seen as an extension of Proposition \ref{prop:stronger}.
\begin{proposition} \label{prop:strongest}
    For $n\ge 4$,
    \[
        B_{n+1}B_{n-1} \ge (B_n)^2 + B_nB_{n-2} + B_{n-1}B_{n-3}.
    \]
\end{proposition}

The proof mostly follows the proof of Proposition \ref{prop:stronger}. Therefore, we will derive Theorem \ref{thm:main} from Proposition \ref{prop:strongest} first. Instead of directly manipulating the inequalities, a proof by induction appears to be more convenient.
\begin{proof}[Proof of Theorem \ref{thm:main}]
We prove by induction. The base cases $n=5,6$ can be checked directly:
\[
B_6B_4=3045\ge 3004=B_5^2+\sum_{i=1}^{5}F_iB_{5-i}^2,
\]
and
\[
B_7B_5=45604\ge 44479=B_6^2+\sum_{i=1}^{6}F_iB_{6-i}^2.
\]

Assuming the conclusion holds for smaller values, we prove it also holds for some $n\ge 7$.
By Proposition 3 and the induction hypothesis applied to
$n-1$ and $n-2$, we have
\begin{align*}
B_{n+1}B_{n-1}
&\ge (B_n)^2 + B_nB_{n-2} + B_{n-1}B_{n-3} \\
&\ge (B_n)^2
    + \sum_{i=0}^{n-1} F_i (B_{n-1-i})^2
    + \sum_{i=0}^{n-2} F_i (B_{n-2-i})^2 \\
&= F_0(B_n)^2
    + F_0(B_{n-1})^2
    + \sum_{j=2}^{n} F_{j-1}(B_{n-j})^2
    + \sum_{j=2}^{n} F_{j-2}(B_{n-j})^2 \\
&= F_0(B_n)^2
    + F_1(B_{n-1})^2
    + \sum_{j=2}^{n} F_j(B_{n-j})^2 \\
&= \sum_{j=0}^{n} F_j(B_{n-j})^2.
\end{align*}
Here, the reindexing is obtained by setting $j=i+1$ in the first sum and $j=i+2$ in the second sum.
The induction step is done, hence the conclusion follows.
\end{proof}

It remains to prove Proposition \ref{prop:strongest}.
\begin{proof}[Proof of Proposition \ref{prop:strongest}]
    Let $R= \{1+s: s\in \mathcal S_{n-3}\}$. Applying Lemma \ref{lem:recurrence} for $t=0,1,2,3,4$ with $t=4$ for
    \[
        B_{n+1}=\sum_{k\in\mathcal S_{n-3}}\left(k^4+4k^3+12k^2+20k+15\right)=\sum_{x\in R}\left(x^4+6x^2+4x+4\right),
    \]
    we rewrite the conclusion as
    \begin{multline*}
        \left(\sum_x x^4+6x^2+4x+4\right) \left(\sum_x x^2+1 \right) \ge \left(\sum_x x^3+3x+1 \right)^2 \\
        + \left(\sum_x x^3+3x+1 \right)\left(\sum_x x\right) + \left(\sum_x x^2+1 \right)\left(\sum_x 1 \right)
    \end{multline*}
where the sums are over $R$.

Proving the inequality is standard and we can follow the same strategy as in Proposition \ref{prop:stronger}, but with more complicated computations.

Let
\[
A(x)=x^4+6x^2+4x+4,\qquad B(x)=x^2+1,\qquad C(x)=x^3+3x+1.
\]
The desired inequality is equivalent to $D\ge 0$, where
\[
D=\left(\sum_x A(x)\right)\left(\sum_x B(x)\right)
-\left(\sum_x C(x)\right)^2
-\left(\sum_x C(x)\right)\left(\sum_x x\right)
-\left(\sum_x B(x)\right)\left(\sum_x 1\right).
\]
Thus
\[
D=\sum_{x,y}\left[A(x)B(y)-C(x)C(y)-C(x)y-B(x)\right].
\]
Doubling and symmetrizing gives
\[
2D=\sum_{x,y}\Big[
A(x)B(y)+A(y)B(x)-2C(x)C(y)
-C(x)y-C(y)x-B(x)-B(y)
\Big].
\]
A direct expansion of the summand gives
\begin{align*}
& A(x)B(y)+A(y)B(x)-2C(x)C(y)-C(x)y-C(y)x-B(x)-B(y) \\
= &
\frac{1}{2}(x-y)^2
\left(2x^2y^2+2x^2+2y^2-10xy-4x-4y+17\right)
+\frac{1}{2}(x+y-4)(4xy+x+y-2).
\end{align*}
It remains to check that both terms are nonnegative.

Since $x,y\ge 2$, we have
\[
x+y-4\ge 0
\qquad\text{and}\qquad
4xy+x+y-2>0.
\]
Also, setting $t=xy\ge 4$,
\[
2x^2y^2+2x^2+2y^2-10xy-4x-4y+17
=(2t^2-10t+17)+(2x^2-4x)+(2y^2-4y).
\]
Now
\[
2x^2-4x\ge 0,\qquad 2y^2-4y\ge 0,
\]
and
\[
2t^2-10t+17
=
2\left(t-\frac{5}{2}\right)^2+\frac{9}{2}>0.
\]
Hence every summand in the expression for $2D$ is nonnegative. Therefore $2D\ge 0$, and so $D\ge 0$. This proves the desired inequality.
\end{proof}

\bibliographystyle{unsrt}
\bibliography{bell}
\end{document}